\documentclass[11pt]{article}
\usepackage{amssymb,amsmath,mathrsfs,txfonts,graphicx,tikz,color}
\usepackage[numbers,sort&compress]{natbib}

\begin{document}
\allowdisplaybreaks
\def\pd#1#2{\frac{\partial#1}{\partial#2}}
\let\oldsection\section
\renewcommand\section{\setcounter{equation}{0}\oldsection}
\renewcommand\thesection{\arabic{section}}
\renewcommand\theequation{\thesection.\arabic{equation}}
\newtheorem{theorem}{\indent Theorem}[section]
\newtheorem{lemma}{\indent Lemma}[section]
\newtheorem{proposition}{\indent Proposition}[section]
\newtheorem{definition}{\indent Definition}[section]
\newtheorem{remark}{\indent Remark}[section]
\newtheorem{corollary}{\indent Corollary}[section]

\title{\LARGE
Critical Sharp Front for Doubly Nonlinear Degenerate Diffusion Equations with Time Delay}
\author{
Tianyuan Xu$^{a}$, Shanming Ji$^{a,}$\thanks{Corresponding author, email:jism@scut.edu.cn},
Ming Mei$^{b,c}$, Jingxue Yin$^d$,
\\
\\
{ \small \it $^a$School of Mathematics, South China University of Technology}
\\
{ \small \it Guangzhou, Guangdong, 510641, P.~R.~China}
\\
{ \small \it $^b$Department of Mathematics, Champlain College Saint-Lambert}
\\
{ \small \it Quebec,  J4P 3P2, Canada, and}
\\
{ \small \it $^c$Department of Mathematics and Statistics, McGill University}
\\
{ \small \it Montreal, Quebec,   H3A 2K6, Canada}
\\
{ \small \it $^d$School of Mathematical Sciences, South China Normal University}
\\
{ \small \it Guangzhou, Guangdong, 510631, P.~R.~China}
}
\date{}

\maketitle

\begin{abstract}
This paper is concerned with the critical sharp traveling wave for doubly nonlinear
diffusion equation with time delay, where the doubly nonlinear degenerate diffusion is defined by $\Big(\big|(u^m)_x\big|^{p-2}(u^m)_x\Big)_x$
with $m>0$ and $p>1$.
The doubly nonlinear diffusion equation is proved to admit a unique
sharp type traveling wave for the degenerate case $m(p-1)>1$, the so-called slow-diffusion case.
This sharp traveling wave associated with the minimal wave speed $c^*(m,p,r)$ is monotonically increasing, where the minimal wave speed satisfies
$c^*(m,p,r)<c^*(m,p,0)$ for any time delay $r>0$.
The sharp front is $C^1$-smooth for $\frac{1}{p-1}<m< \frac{p}{p-1}$,
and piecewise smooth for $m\ge \frac{p}{p-1}$.
Our results indicate
that time delay slows down the minimal traveling wave speed
for the doubly nonlinear degenerate diffusion equations.
The approach adopted for proof is the phase transform method combining the variational method.
The main technical issue for the proof is to overcome
the obstacle caused by the doubly nonlinear degenerate diffusion.

\end{abstract}


\section{Introduction}
This is a continuity of our recent study \cite{Non20} on critical traveling waves for time-delayed degenerate diffusion equation.
Our purpose in the present paper is to study the existence, uniqueness and regularity of the critical sharp traveling wave for the following doubly nonlinear
diffusion equation with time delay
\begin{align}\label{eq-main}
\pd u t=\Big(\big|(u^m)_x\big|^{p-2}(u^m)_x\Big)_x-d(u)+b(u(t-r,x)),\quad x\in \mathbb R,~t>0,
\end{align}
where $p>1$, $m>0$, $u$ is the population density, $b(u(t-r,x))$ is the birth function,
$r\ge0$ is the time delay, and $d(u)$ is the death rate function.
The differential operator $\Big(\big|(u^m)_x\big|^{p-2}(u^m)_x\Big)_x$ is called ``doubly nonlinear''
or non-Newtonian polytropic filtration, see \cite{VazquezJDE17,Jin-JDE} for example.
We focus on the slow diffusion case $m(p-1)>1$
such that sharp type (semi-compactly supported) traveling wave exists and
the initial perturbation propagates at finite speed for the non-delayed case.
The functions $d(s)$ and $b(s)$ satisfy the following conditions:
\begin{enumerate}
\item[(H$_1$)]  Two constant equilibria: $u_-=0$ and $u_+>0$ such that $d(0)=b(0)=0$,  $d(u_+)=b(u_+)$,  $b'(0)>d'(0)\ge 0$, and $d'(u_+)> b'(u_+)\ge0$;
\item[(H$_2$)] Monotonicity: $d(\cdot), \ b(\cdot)\in C^2([0,u_+])$, and  $b'(s)> 0$,  $d'(s)>0$ for $s\in [0,u_+]$.
\end{enumerate}
The assumptions (H$_1$)-(H$_2$) are summarized from
a large number of evolution equations in ecology, such as
the classical Fisher-KPP equation \cite{Fisher};
the well-studied Nicholson's blowflies equation \cite{Gurney80} with
the death function $d_1(u)=\delta u$ or $d_2(u)=\delta u^2$,
the birth function
$$b_1(u)=\tilde{p}u\mathrm{e}^{-au^{\tilde{q}}},
\quad \tilde{p}>0,\quad \tilde{q}>0,\quad a>0;$$
and the Mackey-Glass equation \cite{Mackey} with
the growth function
$$
b_2(u)=\frac{\tilde{p}u}{1+au^{\tilde{q}}}\quad \tilde{p}>0,\quad \tilde{q}>0,\quad a>0.
$$

When $p = 2, \ m=1$, we have the standard heat equation with time delay.
As far as we know, reaction diffusion equations with time delay has
first been studied by Schaaf in \cite{Schaaf}, where he proved the existence of
monotone traveling waves.
The proof was based on sub and super solutions and phase plane techniques.
Since then, the study of traveling wave solutions for reaction diffusion equations with time delay has drawn considerable attention (see, for example, \cite{Mei_LinJDE09,Chern,Faria,Gomez,LLLM} and references therein).
Note that,  the results mentioned above are all for the case that the diffusion term is classical Laplacian.
Choosing $p=2, \ m>1$, we obtain Porous Medium operator, which describes density-dependent dispersal in biological settings.
Here, the important feature of degenerate diffusion equation appears: traveling waves exhibit free boundaries.
In \cite{Non20}, we found the sharp type traveling wave (partially compactly supported) corresponding to the critical wave speed
and obtained the uniqueness of these waves.
Further, we proved that the initial perturbation propagates asymptotically at the same speed
\cite{ArXiv21} and later sharp-oscillatory non-monotone traveling waves was found in \cite{JDE20}.

The sharp type (partially compactly supported) traveling wave solutions are essential in the
analysis of the propagation properties of degenerate diffusion equations.
In many cases, the solutions with (partially) compactly supported initial data
propagate asymptotically at the same speed of the sharp waves,
which also is the minimal admissible traveling wave speed.
This phenomenon was observed by Audrito and V\'azquez \cite{VazquezJDE17}
for doubly nonlinear diffusion equation \eqref{eq-main} without time delay (i.e., $r=0$),
and further the speed was characterized via a variational approach by Benguria and Depassier \cite{Benguria18}.

The existence of traveling waves of \eqref{eq-main} remains to be
technically demanding.
Our main objective is to investigate the structure of the critical sharp waves  and to estimate the corresponding critical speed using the approach of phase transform method with the help of the variational approach  developed  recently in our studies \cite{JDE18,Non20}.
Precisely speaking, we prove that, the doubly nonlinear diffusion equation \eqref{eq-main}  possesses a unique
sharp type traveling wave $\phi(x+c^*t)$ for the degenerate case $m(p-1)>1$, and
such a sharp traveling wave associated with the minimal wave speed $c^*=c^*(m,p,r)$ is monotonically increasing, where the minimal wave speed satisfies
$c^*(m,p,r)<c^*(m,p,0)$ for any time delay $r>0$. Furthermore, we show the optimal regularity of the sharp front $\phi(x+c^*t)$.
That is, when
$\frac{p-1}{m(p-1)-1}$ is integer, then the sharp front $\phi(x+c^*t)$ is $C^{\frac{p-1}{m(p-1)-1}-1}$-smooth with $\phi$ and all its derivatives $\partial^j \phi$ are Lipschitz continuous
for $j=1,\cdots, \frac{p-1}{m(p-1)-1}-1$; while, when $\frac{p-1}{m(p-1)-1}$ is non-integer, then the sharp front $\phi(x+c^*t)$ is
$C^{[\frac{p-1}{m(p-1)-1}]}$-smooth, where $[\frac{p-1}{m(p-1)-1}]$ denotes the largest integer which is less then $\frac{p-1}{m(p-1)-1}$, in particular,
$\phi$ and its all derivatives $\partial^j \phi$ for $j=1,\cdots, [\frac{p-1}{m(p-1)-1}]$ are
$C^{\alpha_{m,p}}$ H\"older continuous with the H\"older exponent $\alpha_{m,p}=\frac{p-1}{m(p-1)-1}-[\frac{p-1}{m(p-1)-1}]$.
This implies that
the sharp front $\phi(x+c^*t)$ is $C^1$-smooth for $\frac{1}{p-1}<m< \frac{p}{p-1}$,
and piecewise smooth for $m\ge \frac{p}{p-1}$.
On the other hand, we also prove
that the time delay $r>0$ slows down the minimal traveling wave speed $c^*=c^*(m,p,r)$
for the doubly nonlinear degenerate diffusion equations.
Finally, let us point out a slightly unexpected phenomenon
related to the doubly nonlinear operator.
The main difficulty lies in the asymptotic behavior of the phase function $\tilde\psi(\phi)$
defined for the sharp type traveling wave $\phi(\xi)$ by regarding $\psi(\xi):=|(\phi^m(\xi))'|^{p-2}(\phi^m(\xi))'$
as a function of $\phi$.
Its asymptotic behavior near the positive equilibrium $u_+$ for the degenerate case $p\in(1,2)$ is quite different from the
case $p=2$.

The paper is organized as follows. In Section 2, we state our main results.
We defer to Section 3 all the detailed proofs. Section 4 is the brief derivation of models we treat.

\section{Main results}
We consider the doubly nonlinear degenerate diffusion equation with time delay \eqref{eq-main}.
We are looking for the traveling wave solutions of sharp type that
connect the two equilibria $u_-=0$ and $u_+=:K$.
Under the hypotheses (H$_1$)-(H$_2$),
the birth function $b(u)$ is monotonically increasing on $[u_-,u_+]=:[0,K]$.
Let $\phi(\xi)$, where $\xi=x+ct$ and $c>0$,
be the traveling wave solution of  \eqref{eq-main},
we get (we write $\xi$ as $t$ for the sake of simplicity)
\begin{align}\label{eq-tw}
\begin{cases}
\displaystyle
c\phi'(t)=(|(\phi^m)'(t)|^{p-2}(\phi^m)'(t))'-d(\phi(t))+b(\phi(t-cr)),\quad t\in\mathbb R,\\
\phi(-\infty)=0, \quad \phi(+\infty)=K.
\end{cases}
\end{align}

Since \eqref{eq-tw} has singularity or degeneracy, we employ the following definition of
sharp and smooth traveling waves.
Here are some notations used throughout this paper:
$$
C_\mathrm{unif}^\mathrm{b}(\mathbb R):=
\{\phi\in C(\mathbb R)\cap L^\infty(\mathbb R);\phi \text{~is uniformly continuous on~}\mathbb R\},
$$
and
$$W_\mathrm{loc}^{1,p}(\mathbb R):=\{\phi; \phi\in W^{1,p}(\Omega)
\text{~for any compact subset~}\Omega\subset\mathbb R\}.$$

\begin{definition} \label{de-semi}
A profile function $\phi(t)$ is said to be a traveling wave solution
of \eqref{eq-tw} if $\phi\in C_\mathrm{unif}^\mathrm{b}(\mathbb R)$, $0\le \phi(t)\le K:=u_+$,
$\phi(-\infty)=0$, $\phi(+\infty)=K$,
$\phi^m\in W_{\mathrm{loc}}^{1,p}(\mathbb R)$,
$\phi(t)$ satisfies \eqref{eq-tw} in the sense of distributions.
The traveling wave $\phi(t)$ is said to be of sharp type if
the support of $\phi(t)$ is semi-compact, i.e.,
$\text{supp}\,\phi=[t_0,+\infty)$ for some $t_0\in\mathbb R$,
$\phi(t)>0$ for $t>t_0$.
On the contrary, the traveling wave $\phi(t)$
is said to be of smooth type if $\phi(t)>0$ for all $t\in\mathbb R$.
\end{definition}

Without loss of generality, we may always shift $t_0$ to $0$
for the sharp type traveling wave.
Therefore, a sharp type traveling wave $\phi(t)$ is
a special solution such that $\phi(t)=0$ for $t\le0$, and $\phi(t)>0$ for $t>0$.

For any given $m>0$, $p>1$, such that $m(p-1)>1$, and $r\ge0$, we define the critical (or minimal) wave speed
$c^*(m,p,r)$ for the degenerate diffusion equation \eqref{eq-tw} as follows
\begin{equation} \label{eq-def}
c^*(m,p,r)
:=\inf\{c>0; \eqref{eq-tw} \text{~admits increasing traveling waves with speed $c$}\}.
\end{equation}
For the case without time delay and with degenerate diffusion (i.e. $m(p-1)>1$ and $r=0$),
it is proved by Benguria and Depassier in \cite{Benguria18} that
\begin{equation} \label{eq-cstar0}
c^*(m,p,0)=\sup_{g\in \mathscr{D}}
\int_0^K \frac{p}{(p-1)^{(p-1)/p}}(-g'(\phi))^\frac{1}{p}(g(\phi))^\frac{p-1}{p}
(m\phi^{m-1}(b(\phi)-d(\phi)))^\frac{p-1}{p}\mathrm{d}\phi,
\end{equation}
where $\mathscr{D}=\{g\in C^1([0,K]);\int_0^K g(s)\mathrm{d}s=1,g(s)>0,g'(s)<0,\forall s\in(0,K)\}$.

In this paper we show that \eqref{eq-tw} admits a unique
sharp type traveling wave, and the sharp traveling wave is monotonically increasing
and corresponding to the minimal wave speed $c^*(m,p,r)$,
and further $c^*(m,p,r)<c^*(m,p,0)$ for any time delay $r>0$.
As a consequence, the time delay slows down the minimal traveling wave speed
for the doubly nonlinear degenerate diffusion equations.

Our main results are as follows.

\begin{theorem}[Critical Sharp Traveling Wave] \label{th-existence}
Assume that $d(s)$ and $b(s)$ satisfy (H$_1$)-(H$_2$), and $m>0$, $p>1$, $r\ge0$, such that $m(p-1)>1$.
There exists a unique $c^*=c^*(m,p,r)>0$ defined in \eqref{eq-def}
satisfying $c^*(m,p,r)<c^*(m,p,0)$ for any time delay $r>0$,
such that \eqref{eq-tw} admits a unique (up to shift)
sharp traveling wave $\phi(x+c^*t)$ with speed $c^*$,
which is the critical traveling wave of \eqref{eq-tw} and is monotonically increasing.
Moreover, any other traveling wave solution must be smooth and correspond to speed $c>c^*(m,p,r)$.
\end{theorem}

\begin{theorem}[Regularity of Sharp Wave] \label{th-sharp}
Assume that the conditions in Theorem \ref{th-existence} hold.
Let $\gamma_{m,p}$ be the largest integer that is smaller than $\frac{p-1}{m(p-1)-1}$, i.e.,
$$
\gamma_{m,p}:=
\begin{cases}
\frac{p-1}{m(p-1)-1}-1, \quad& \text{if~} \frac{p-1}{m(p-1)-1} \text{is an integer},\\
[\frac{p-1}{m(p-1)-1}], \quad& \text{if~} \frac{p-1}{m(p-1)-1} \text{is not an integer},
\end{cases}
$$
and denote $\alpha_{m,p}:=\frac{p-1}{m(p-1)-1}-\gamma_{m,p}\in(0,1]$.
Then the optimal regularity of sharp wave $\phi(\xi)$ is
$\phi\in C^{\gamma_{m,p},\alpha_{m,p}}(\overline{\mathbb R})$,
where $C^{\gamma_{m,p},\alpha_{m,p}}(\overline{\mathbb R})$ is the function space defined as:
if $\frac{p-1}{m(p-1)-1}$ is an integer, then $\alpha_{m,p}=1$, and
\begin{eqnarray}
C^{\gamma_{m,p},1}(\overline{\mathbb R}):=\Big\{
\phi \in C^{\gamma_{m,p}}(\overline{\mathbb R})&\Big| & \partial^j \phi, \mbox{ for } j=0,1,\cdots, \frac{p-1}{m(p-1)-1}-1, \notag \\
& &\mbox{ are Lipschitz continuous}
\Big\};
\end{eqnarray}
while if $\frac{p-1}{m(p-1)-1}$ is not an integer, then $0<\alpha_{m,p}<1$, and
\begin{eqnarray}
C^{\gamma_{m,p},\alpha_{m,p}}(\overline{\mathbb R}):=\Big\{
\phi \in C^{\gamma_{m,p}}(\overline{\mathbb R})&\Big| & \partial^j \phi, \mbox{ for } j=0,1,\cdots, \Big[\frac{p-1}{m(p-1)-1}\Big], \notag \\
& & \mbox{ are } C^{\alpha_{m,p}} \mbox{  H\"older continuous}
\Big\}.
\end{eqnarray}
\end{theorem}

\begin{figure}[htb]
\begin{center}
\begin{tikzpicture}[scale=0.8,domain=-6:6]
\def\axl{-2} \def\axr{5} \def\ay{4} \def\pa{\ay/180} \def\pb{3} \def\pc{2}
\def\qa{-0.1531} \def\qb{1.1531} \def\ra{1.3606} \def\rb{1.2789}
\draw[line width=.6pt] ({\axl},0)--(\axr,0);
\draw[line width=.6pt] ({\axl},\ay)--(\axr,\ay);
\draw[line width=.6pt] ({\axl+8},0)--(\axr+8,0);
\draw[line width=.6pt] ({\axl+8},\ay)--(\axr+8,\ay);
\draw[color=blue,line width=.7pt]
    plot[domain=2+8:5+8, samples=144, smooth] (\x,{\ay/2+\pa*atan(\pb*(\x-8-2))});
\draw[color=blue,line width=.7pt]
    plot[domain=0+8:1+8, samples=144, smooth] (\x,{(\x-8)*(\x-8)/3});
\draw[color=blue,line width=.7pt]
    plot[domain=1+8:2+8, samples=144, smooth]
    (\x,{1/3+2/3*(\x-8-1)+\qa*(\x-8-1)*(\x-8-1)+\qb*(\x-8-1)^3});
\draw[color=red,line width=.7pt]
    plot[domain=1:5, samples=144, smooth] (\x,{\ay/2+\pa*atan(\pb*(\x-1))});
\draw[color=red,line width=.7pt]
    plot[domain=0:0.5, samples=144, smooth] (\x,{\x*3/2});
\draw[color=red,line width=.7pt]
    plot[domain=0.5:1, samples=144, smooth]
    (\x,{3/4+3/2*(\x-1/2)+\ra*(\x-1/2)*(\x-1/2)+\rb*(\x-1/2)^3});
\draw[color=red,line width=.7pt] ({\axl},0)--(0,0);
\draw[color=blue,line width=.7pt] ({\axl+8},0)--(0+8,0);
\node at (2.8,1.4) {($a$) $m\ge\frac{p}{p-1}$};
\node at (8.2,2.4) {($b$) $m\in(\frac{1}{p-1},\frac{p}{p-1})$};
\end{tikzpicture}
\end{center}
\caption{Traveling waves: ($a$) non-$C^1$ sharp type for $m\ge\frac{p}{p-1}$;
($b$) $C^1$ sharp type for $m\in(\frac{1}{p-1},\frac{p}{p-1})$.}
\label{fig-1}
\end{figure}
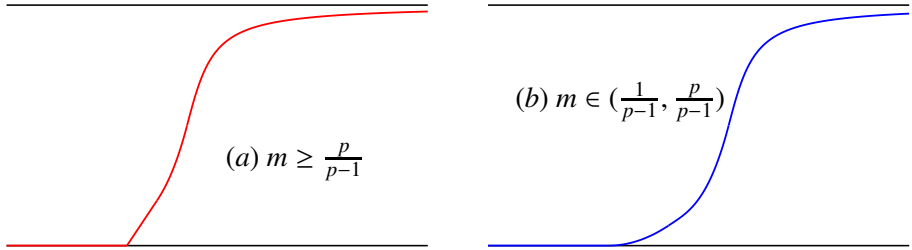

\begin{remark}
If $m\ge\frac{p}{p-1}$, then the sharp traveling wave is not $C^1$ smooth;
while if $m\in(\frac{1}{p-1},\frac{p}{p-1})$, then the sharp traveling wave is $C^1$ smooth. See Figure \ref{fig-1}.
\end{remark}

\section{Proof of the main results}
For any given $m>0$, $p>1$, and $r>0$, such that $m(p-1)>1$,
we solve \eqref{eq-tw} locally for any $c>0$
and then we single out a special one that is a sharp traveling wave with critical wave speed.
First, noticing that the sharp wave solution $\phi(t)=0$ for $t\le0$
and then $\phi(t-cr)=0$ for $t\in[0,cr)$, \eqref{eq-tw} is locally reduced to
\begin{equation} \label{eq-semi-1}
\begin{cases}
c\phi'(t)=(|(\phi^m(t))'|^{p-2}(\phi^m(t))')'-d(\phi(t)), \\
\phi(0)=0, \quad (\phi^m)'(0)=0, \quad t\in(0,cr),
\end{cases}
\end{equation}
whose solutions are not unique and we choose the maximal one
such that $\phi(t)>0$ for $t\in(0,cr)$ as shown in the following lemma.
Here, $(\phi^m)'(0)=0$ is a necessary and sufficient condition
such that the zero extension of $\phi(t)$ to the left
satisfies \eqref{eq-tw} locally near $0$ in the sense of distributions.

The proof follows from the similar outline as in \cite{Non20},
the difference lies in the asymptotic behavior of $\psi(t):=|(\phi^m(t))'|^{p-2}(\phi^m(t))'$ in the singular phase plane
of $(\phi,\psi)$ for the sharp wave solution $\phi(t)$.
Here we mainly sketch the proofs that have differences for the sake of simplicity.

\begin{lemma} \label{le-semi-1}
For any $c>0$, the degenerate ODE \eqref{eq-semi-1} admits a unique maximal solution
$\phi_c^1(t)$ on $(0,cr)$ such that $\phi_c^1(t)>0$ on $(0,cr)$ and
\begin{equation} \label{eq-expansion}
\phi(t)=\Big(\frac{c^\frac{1}{p-1}(m(p-1)-1)}{m(p-1)}\Big)^\frac{p-1}{m(p-1)-1}\cdot t_+^\frac{p-1}{m(p-1)-1}
+o(t_+^\frac{p-1}{m(p-1)-1}), \quad \text{as~} \ t\to0^+.
\end{equation}
\end{lemma}
{\it\bfseries Proof.}
A positive function $\phi(t)>0$ on $(0,cr)$ is a solution to the degenerate ODE \eqref{eq-semi-1}
satisfies the following singular differential system on $(0,cr)$
\begin{equation} \label{eq-zsys}
\begin{cases}
\displaystyle
\phi'(t)=\frac{\psi^\frac{1}{p-1}(t)}{m\phi^{m-1}(t)},
\\ \displaystyle
\psi'(t)=c\frac{\psi^\frac{1}{p-1}(t)}{m\phi^{m-1}(t)}+d(\phi(t)),
\end{cases}
\end{equation}
with $\psi(t):=|(\phi^m(t))'|^{p-2}(\phi^m(t))'$.
We seek for a solution to \eqref{eq-zsys} such that $\phi(t)>0$ and $\psi(t)>0$ for $t\in(0,cr)$,
with $\psi(0)=0$ and $\phi(0)=0$.
The system \eqref{eq-zsys} has singularity at
some points where $\phi(t)=0$.
Therefore, we solve \eqref{eq-zsys} with the initial condition
$(\phi_\varepsilon(0),\psi_\varepsilon(0))=(\varepsilon^2,\varepsilon)$,
whose local existence and uniqueness are certain according the classical phase plane analysis method.

For any trajectory $(\phi(t),\psi(t))$ such that $\phi(t)>0$ and $\psi(t)>0$ in some interval,
we can make change of variables such that we
take $\phi$ as an independent variable and regard $\psi$ as a function of $\phi$,
denoted by $\tilde \psi(\phi)$,
in the phase plane of $(\phi,\psi)$
since $\phi'(t)=\frac{\psi^{1/(p-1)}(t)}{m\phi^{m-1}(t)}>0$.
The function corresponding to $(\phi_\varepsilon,\psi_\varepsilon)$ is denoted by $\tilde \psi_\varepsilon(\phi)$.
Comparison principle or the analysis of the trajectories shows that
$\tilde \psi_\varepsilon(\phi)$ is monotone increasing with respect to $\varepsilon>0$.
The limiting function $(\phi(t),\psi(t))$ as $\varepsilon\to0^+$ is the
maximal solution to \eqref{eq-semi-1}.

Asymptotic analysis shows that (note that $m(p-1)>1$)
$$
\tilde \psi(\phi)=c\phi+o(\phi), \quad \text{as~} \ \phi\to0^+,
$$
or equivalently,
\begin{equation} \label{eq-ztildepsi}
\psi(t)=|(\phi^m(t))'|^{p-2}(\phi^m(t))'=c\phi(t)+o(\phi(t)), \quad \text{as~} \ t\to0^+.
\end{equation}
Furthermore, $\phi^m(t)>0$ on $(0,cr)$ with $\phi^m(0)=0$ is the maximal solution to the following
singular first order differential equation
$$
(\phi^m(t))'=c^\frac{1}{p-1}(\phi^m(t))^\frac{1}{m(p-1)}+o((\phi^m(t))^\frac{1}{m(p-1)}), \quad \text{as~} \ t\to0^+.
$$
Therefore,
$$
\phi^m(t)=\Big(\frac{c^\frac{1}{p-1}(m(p-1)-1)}{m(p-1)}\Big)^\frac{m(p-1)}{m(p-1)-1}\cdot t_+^\frac{m(p-1)}{m(p-1)-1}
+o(t_+^\frac{m(p-1)}{m(p-1)-1}), \quad \text{as~} \ t\to0^+.
$$
The proof of \eqref{eq-expansion} is completed.
$\hfill\Box$

Next, let $\phi_c^2(t)$ be the solution of the following initial value second order ODE problem
\begin{equation} \label{eq-semi-2}
\begin{cases}
c\phi'(t)=(|(\phi^m(t))'|^{p-2}(\phi^m(t))')'-d(\phi(t))+b(\phi_c^1(t-cr)), \quad t\in(cr,2cr),\\
\phi(cr)=\phi_c^1(cr), \quad
\phi'(cr)=(\phi_c^1)'(cr).
\end{cases}
\end{equation}
The problem \eqref{eq-semi-2} is locally solvable and has no singularity near $t=cr$
since $\phi_c^1(cr)>0$.
The above steps can be continued unless $\phi_c^k(t)$ blows up
or decays to zero in finite time for some $k\in \mathbb N^+$.
Let $\phi_c(t)$ be the connecting function of those functions on each step, i.e.,
\begin{equation} \label{eq-semi}
\phi_c(t)=
\begin{cases}
\phi_c^1(t), \quad &t\in[0,cr),\\
\phi_c^2(t), \quad &t\in[cr,2cr),\\
\vdots\\
\phi_c^k(t), \quad &t\in[(k-1)cr,kcr),\\
\vdots
\end{cases}
\end{equation}
for some finite steps such that $\phi_c(t)$ blows up or decays to zero,
or for infinite steps such that $\phi_c(t)$ is defined on $(0,+\infty)$
and zero extended to $(-\infty,0)$ for convenience.

\begin{lemma} \label{le-semi-decay}
For any given $r>0$, there exist two numbers $\overline c>\underline c>0$ such that
if $0<c\le \underline c$, then
$\phi_c(t)$ decays to zero in finite time;
if $c\ge \overline c$, then
$\phi_c(t)$ grows up to $+\infty$ as $t$ tends to $+\infty$.
\end{lemma}
{\it\bfseries Proof.}
The above assertions for the special case of $m>1$ and $p=2$ are proved in \cite{JDE20}.
Here we consider the generalized phase plane corresponding to the following dynamical system
\begin{equation} \label{eq-zphase}
\begin{cases}
\displaystyle
\phi'(t)=\frac{\psi^\frac{1}{p-1}(t)}{m\phi^{m-1}(t)},
\\ \displaystyle
\psi'(t)=c\frac{\psi^\frac{1}{p-1}(t)}{m\phi^{m-1}(t)}+d(\phi(t))-b(\phi(t-cr)),
\end{cases}
\end{equation}
with $\psi(t):=|(\phi^m(t))'|^{p-2}(\phi^m(t))'$.
For the local solution $\phi_c(t)$ on its strictly monotone increasing subinterval,
we take $\phi_c$ as an independent variable and regard $\psi_c$ as a function of $\phi_c$, denoted by
\begin{equation} \label{eq-psi}
\tilde \psi_c(\phi_c)=\psi_c(t^{-1}(\phi_c)), \text{~such~that~} t^{-1}(\phi_c) \text{~is~the~inverse~function~of~} \phi_c(t).
\end{equation}
We may drop the subscripts in $\phi_c$, $\psi_c$, $\tilde\psi_c$, and simply write $\phi$, $\psi$, $\tilde\psi$,
for a given $c>0$.
Define
\begin{equation} \label{eq-phicr}
\phi_{cr}:=\inf_{\theta\in[0,\phi]}\Big\{\int_\theta^\phi\frac{ms^{m-1}}{\tilde\psi^\frac{1}{p-1}(s)}\mathrm{d}s\le cr\Big\}.
\end{equation}
Then the function $\tilde\psi(\phi)$ satisfies
\begin{equation} \label{eq-zphasepsi}
\begin{cases}
\displaystyle
\frac{\mathrm{d}\tilde\psi}{\mathrm{d}\phi}=c-\frac{m\phi^{m-1}\cdot(b(\phi_{cr})-d(\phi))}{\tilde\psi^\frac{1}{p-1}(\phi)},\\
\tilde\psi(0)=0, \quad \tilde\psi(\phi)>0 \text{~for~} \phi\in(0,\phi^*),
\end{cases}
\end{equation}
where $\phi^*=\phi_c(t^*)$ and $(0,t^*)$ is the maximum interval such that $\phi_c(t)$ is strictly monotone increasing.
The rest of the proof follows similarly as the proofs of Lemma 3.2 and Lemma 3.3 in \cite{JDE20}.
$\hfill\Box$

The local solution $\phi_c(t)$ may grow beyond the positive equilibrium $K>0$ or
decay to zero in finite interval.
The sharp traveling wave is the special one (the uniqueness will be proved) such that
$\phi_c(t)$ exists globally and is monotone increasing on $(0,+\infty)$,
together with the speed $c$ being identical to the critical wave speed $c^*(m,p,r)$.
The existence and other properties of the sharp wave
for the case of $m>1$ and $p=2$ are proved in \cite{JDE20,Non20}.
Specifically, the existence of sharp wave
in the above settings is prove in \cite{JDE20} through a continuous argument
for general non-monotone birth function $b(s)$;
and the uniqueness is prove in \cite{Non20} via monotone dependence of
$\phi_c(t)$ with respect to $c$ for monotone birth function.

\begin{lemma} \label{le-dependent}
The solution $\phi_c(t)$ is locally continuously dependent on $c$
and is strictly monotonically increasing with respect to $c$ on their joint existence interval.
\end{lemma}
{\it\bfseries Proof.}
The proof is similar to that in Lemma 3.4 and Lemma 3.6 in \cite{Non20}.
Here we omit the details.
$\hfill\Box$

\begin{lemma} \label{le-cstar}
There exists a unique number $c^*=c^*(m,p,r)>0$ such that
$\phi_{c^*}(t)$ is strictly increasing on $(0,+\infty)$ with
$\phi_{c^*}(+\infty)=K$,  and the function $\phi_{c^*}(t)$ is the unique traveling wave solution of sharp type.
The speed of any smooth traveling wave is greater than $c^*(m,p,r)$,
and no traveling waves $\phi(x+ct)$ exist when $c\le c^*$. Namely,
$c^*$ is the minimal admissible traveling wave speed.
\end{lemma}
{\it\bfseries Proof.}
This is proved in a similar way as Lemma 3.7 and Lemma 3.9 in \cite{Non20}.
Here we need to note the following asymptotic behavior near zero of the phase function $\tilde \psi(\phi)$,
defined as \eqref{eq-psi}
for any traveling wave solution $\phi(t)$ with $\psi(t):=|(\phi^m(t))'|^{p-2}(\phi^m(t))'$:

(i) if $\phi(t)$ is a sharp traveling wave with speed $c$,
then $\tilde\psi(\phi)=c\phi+o(\phi)$, as $\phi\to0^+$, according to Lemma \ref{le-semi-1};

(ii) if $\phi(t)$ is a smooth traveling wave with speed $c$,
then
$$
\tilde\psi(\phi)=\Big(\frac{m(b'(0)e^{-\lambda cr}-d'(0))}{c}\Big)^{p-1}\phi^{m(p-1)},
\quad \text{as~} \ \phi\to0^+,
$$
where $\lambda>0$ is the unique root of the equation $c\lambda+d'(0)=b'(0)e^{-\lambda cr}$.

Suppose that $\hat\phi(t)$ is a smooth traveling wave with speed $c>0$.
The phase function corresponds to $\hat\phi(t)$ is denoted by $\hat\psi(\phi)$.
Let $\phi_c(t)$ be the local solution of sharp type with the same speed
and $\tilde\psi(\phi)$ be its phase function.
Locally near zero,
$$
\hat\psi(\phi)\sim\Big(\frac{m(b'(0)e^{-\lambda cr}-d'(0))}{c}\Big)^{p-1}\phi^{m(p-1)}
<c\phi\sim \tilde\psi(\phi), \quad \text{as~} \ \phi\to0^+,
$$
since $m(p-1)>1$.
Therefore, the monotone dependence (proved in a similar way as Lemma 3.6 in \cite{Non20})
shows that $\tilde\psi(\phi)>\hat\psi(\phi)$ for $\phi\in(0,K]$.
That is, the local solution $\phi_c(t)$ grows beyond $K$ in finite interval
and there holds $c>c^*(m,p,r)$.
$\hfill\Box$

\medskip

We also need to describe the asymptotic behavior of the phase function $\tilde\psi(\phi)$ near the positive equilibrium $K$,
which is important to the variational characterization of critical wave speed.

\begin{lemma} \label{le-asymp}
The phase function $\tilde \psi(\phi)$ defined in \eqref{eq-psi}
for the unique sharp traveling wave solution $\phi(t)$ satisfies the following asymptotic expansion:

(i) if $p=2$, then
$$
\tilde\psi(\phi)=\kappa_2(K-\phi)+o(|K-\phi|), \quad \text{as~} \ \phi\to K^-,
$$
where $\lambda>0$ is the unique positive root of
$$mK^{m-1}\lambda^2+c\lambda+b'(K)e^{\lambda cr}-d'(K)=0,$$
and $\kappa_2:=mK^{m-1}\lambda$;

(ii) if $p>2$, then
$$
\tilde\psi(\phi)=\kappa_p(K-\phi)^{p-1}+o(|K-\phi|^{p-1}), \quad \text{as~} \ \phi\to K^-,
$$
where $\lambda>0$ is the unique positive root of
$$c\lambda+b'(K)e^{\lambda cr}-d'(K)=0,$$
and $\kappa_p:=(mK^{m-1}\lambda)^{p-1}$ for $p>2$;

(iii) if $1<p<2$, then
$$
\tilde\psi(\phi)=\kappa_p(K-\phi)^\frac{2(p-1)}{p}+o(|K-\phi|^\frac{2(p-1)}{p}), \quad \text{as~} \ \phi\to K^-,
$$
where $\lambda>0$ is the unique positive root of
$$
m^{p-1}K^{(m-1)(p-1)}\frac{2(p-1)}{p}\Big(\frac{p}{2-p}\lambda\Big)^{p}+b'(K)-d'(K)=0.
$$
and $\kappa_p:=m^{p-1}K^{(m-1)(p-1)}\Big(\frac{p}{2-p}\lambda\Big)^{p-1}$ for $p\in(1,2)$.
\end{lemma}
{\it\bfseries Proof.}
For $p\ge2$, we utilize the ansatz of expansion $\tilde\psi(\phi)\sim \kappa(K-\phi)^{p-1}$
and $\phi(t)-K\sim -\mu \mathrm{e}^{-\lambda t}$ as $t\to+\infty$ and $\phi\to K^-$
for some positive constants $\kappa$, $\mu$ and $\lambda$.
Further,
$d(\phi)-d(K)\sim d'(K)(\phi-K)$,
and $b(\phi_{cr})-b(K)\sim b'(K)(\phi_{cr}-K)$, as $t\to+\infty$,
where $\phi_{cr}=\phi(\cdot-cr)$.
Noticing that $b(K)=d(K)$, we see that
\begin{align*}
b(\phi_{cr})-d(\phi)&\sim b'(K)(\phi_{cr}-K)-d'(K)(\phi-K) \\
&\sim (b'(K)-d'(K))(\phi-K)+b'(K)(\phi_{cr}-\phi) \\
&\sim (b'(K)-d'(K))(\phi-K)+b'(K)(\mathrm{e}^{\lambda cr}-1)(\phi-K) \\
&\sim (b'(K)\mathrm{e}^{\lambda cr}-d'(K))(\phi-K),
\end{align*}
as $t\to+\infty$ since $\frac{\phi_{cr}-\phi}{\phi-K}\sim \mathrm{e}^{\lambda cr}-1$.
According to \eqref{eq-zphase} and \eqref{eq-zphasepsi}, near $K$, $\tilde\psi$ behaves similar as
\begin{equation} \label{eq-zphasepsi-K}
\begin{cases}
\displaystyle
\frac{\mathrm{d}\tilde\psi}{\mathrm{d}\phi}\sim
c-\frac{mK^{m-1}\cdot(b'(K)\mathrm{e}^{\lambda cr}-d'(K))(\phi-K)}{\tilde\psi^\frac{1}{p-1}(\phi)},\\
\tilde\psi(K)=0, \quad \tilde\psi(\phi)>0 \text{~for~} \phi\in(0,K).
\end{cases}
\end{equation}

For the special case $p=2$, $\frac{\mathrm{d}\tilde\psi}{\mathrm{d}\phi}\sim -\kappa$, and
the singular ODE \eqref{eq-zphasepsi-K} admits a solution satisfying the expansion $\tilde\psi(\phi)\sim \kappa(K-\phi)$
provided that $\kappa>0$ is the unique positive root of the following equation
\begin{equation} \label{eq-zkappa}
-\kappa=c+\frac{mK^{m-1}\cdot(b'(K)\mathrm{e}^{\lambda cr}-d'(K))}{\kappa}.
\end{equation}
Additionally, a necessary condition for the characteristic value $\lambda>0$ of
the traveling wave of \eqref{eq-tw} satisfying $\phi(t)-K\sim -\mu \mathrm{e}^{-\lambda t}$ as $t\to+\infty$ is
\begin{equation} \label{eq-zlambda}
mK^{m-1}\lambda^2+c\lambda+b'(K)\mathrm{e}^{\lambda cr}-d'(K)=0.
\end{equation}
Moreover, according to the asymptotic expansions $\tilde\psi(\phi)\sim \kappa(K-\phi)$
and $\phi(t)-K\sim -\mu \mathrm{e}^{-\lambda t}$,
noticing that $\tilde\psi(\phi)=\psi(t)=(\phi^m(t))'\sim mK^{m-1}\phi'(t)$ for $p=2$,
there must hold
\begin{equation} \label{eq-zkappa2}
\kappa=mK^{m-1}\lambda.
\end{equation}
Since $d'(K)>b'(K)\ge0$, the characteristic equation \eqref{eq-zlambda} admits a unique positive root $\lambda>0$,
and then \eqref{eq-zkappa} is equivalent to \eqref{eq-zkappa2}.

For the case $p>2$, $\frac{\mathrm{d}\tilde\psi}{\mathrm{d}\phi}=o(1)$ as $\phi\to K^-$,
and in this situation \eqref{eq-zkappa} reads as
\begin{equation} \label{eq-zkappa-p}
0=c+\frac{mK^{m-1}\cdot(b'(K)\mathrm{e}^{\lambda cr}-d'(K))}{\kappa^\frac{1}{p-1}}.
\end{equation}
The characteristic equation \eqref{eq-zlambda} now is
\begin{equation} \label{eq-zlambda-p}
c\lambda+b'(K)\mathrm{e}^{\lambda cr}-d'(K)=0.
\end{equation}
Similar to \eqref{eq-zkappa2} the relation between the expansions
$\tilde\psi(\phi)\sim \kappa(K-\phi)^{p-1}$ and $\phi(t)-K\sim -\mu \mathrm{e}^{-\lambda t}$
implies
\begin{align*}
\tilde\psi(\phi)=\psi(t)=|(\phi^m(t))'|^{p-2}(\phi^m(t))'\sim (mK^{m-1}\phi'(t))^{p-1}
\sim \kappa(K-\phi)^{p-1}.
\end{align*}
That is,
\begin{equation} \label{eq-zkappa2-p}
\kappa=(mK^{m-1}\lambda)^{p-1}.
\end{equation}
The characteristic equation \eqref{eq-zlambda-p} has a unique positive root $\lambda>0$
and \eqref{eq-zkappa-p} is equivalent to \eqref{eq-zkappa2-p} in this case.

The case of $1<p<2$ is quite different, we utilize the ansatz $\tilde\psi(\phi)\sim \kappa(K-\phi)^\frac{2(p-1)}{p}$
and $\phi(t)-K\sim -\mu (1+\lambda t)^{-\frac{p}{2-p}}$ as $t\to+\infty$ and $\phi\to K^-$
for some positive constants $\kappa$, $\mu$ and $\lambda$.
It should be addressed that the sharp traveling wave approached the positive equilibrium $K$ algebraically
instead of exponentially.
Note that $\frac{2(p-1)}{p}\in(0,1)$, $\frac{2(p-1)}{p}-1=1-\frac{2}{p}$, and $\frac{p}{2-p}\in(1,+\infty)$ for $p\in(1,2)$.
The algebraical decay behaves differently from the exponential decay, such that
\begin{align*}
b(\phi_{cr})-d(\phi)&\sim b'(K)(\phi_{cr}-K)-d'(K)(\phi-K) \\
&\sim (b'(K)-d'(K))(\phi-K)+b'(K)(\phi_{cr}-\phi) \\
&\sim (b'(K)-d'(K))(\phi-K),
\end{align*}
as $t\to+\infty$ since $\frac{\phi_{cr}-\phi}{\phi-K}=o(1)$.
The singular ODE \eqref{eq-zphasepsi-K} now shows that
\begin{equation} \label{eq-zkappa-p2}
\frac{2(p-1)}{p}\cdot \kappa=-\frac{mK^{m-1}\cdot(b'(K)-d'(K))}{\kappa^\frac{1}{p-1}}.
\end{equation}
For the degenerate case such that
$\phi(t)-K\sim -\mu (1+\lambda t)^{-\frac{p}{2-p}}$ (assume $\mu=1$ by rescaling), we have
\begin{align*}
(|(\phi^m&(t))'|^{p-2}(\phi^m(t))')'=
(m^{p-1}\phi^{(m-1)(p-1)}|\phi'|^{p-2}\phi')' \\
&=m^{p-1}\phi^{(m-1)(p-1)}(p-1)|\phi'|^{p-2}\phi''+m^{p-1}(m-1)(p-1)\phi^{(m-1)(p-1)-1}|\phi'|^{p} \\
&\sim m^{p-1}\phi^{(m-1)(p-1)}(p-1)|\phi'|^{p-2}\phi'',
\end{align*}
since $|\phi'|^{p}=o(|\phi'|^{p-2}\phi'')$ according to $\frac{p}{2-p}>1$.
Moreover, $\phi'=o(|K-\phi|)=o(b(\phi_{cr})-d(\phi))$, and then the characteristic equation of \eqref{eq-tw} is
$$
-m^{p-1}\phi^{(m-1)(p-1)}(p-1)|\phi'|^{p-2}\phi''\sim b(\phi_{cr})-d(\phi)\sim (b'(K)-d'(K))(\phi-K),
$$
which means
$$
m^{p-1}K^{(m-1)(p-1)}(p-1)\Big(\frac{p}{2-p}\lambda\Big)^{p-2}\frac{p}{2-p}\frac{2}{2-p}\lambda^2+b'(K)-d'(K)=0.
$$
That is,
\begin{equation} \label{eq-zlambda-p2}
m^{p-1}K^{(m-1)(p-1)}\frac{2(p-1)}{p}\Big(\frac{p}{2-p}\lambda\Big)^{p}+b'(K)-d'(K)=0.
\end{equation}
Moreover, according to the asymptotic expansions $\tilde\psi(\phi)\sim \kappa(K-\phi)^\frac{2(p-1)}{p}$
and $\phi(t)-K\sim -(1+\lambda t)^{-\frac{p}{2-p}}$, we have
\begin{align*}
\tilde\psi(\phi)=&|(\phi^m(t))'|^{p-2}(\phi^m(t))'\sim (mK^{m-1}\phi'(t))^{p-1} \\
\sim& m^{p-1}K^{(m-1)(p-1)}\Big(\frac{p}{2-p}(1+\lambda t)^{-\frac{2}{2-p}}\lambda\Big)^{p-1} \\
\sim& \kappa((1+\lambda t)^{-\frac{p}{2-p}})^\frac{2(p-1)}{p}.
\end{align*}
Therefore,
\begin{equation} \label{eq-zkappa2-p2}
\kappa=m^{p-1}K^{(m-1)(p-1)}\Big(\frac{p}{2-p}\lambda\Big)^{p-1}.
\end{equation}
The characteristic equation \eqref{eq-zlambda-p2} has a unique positive root $\lambda$,
and then \eqref{eq-zkappa2-p2} is identical to \eqref{eq-zkappa-p2}.
$\hfill\Box$

\medskip

The critical wave speed $c^*(m,p,0)$ for non-delayed case is characterized via a variational approach
by Benguria and Depassier \cite{Benguria-Variational,Benguria,Benguria18}.
For time-delayed case, we utilize the variational characterization method
to show the dependence of the critical wave speed $c^*(m,p,r)$ with respect to the time delay $r$.

\begin{lemma} \label{le-cstarr}
The minimal traveling wave speed $c^*(m,p,r)$ for the time delay $r>0$
is strictly smaller than that without time delay,
i.e., $c^*(m,p,r)<c^*(m,p,0)$.
\end{lemma}
{\it\bfseries Proof.}
Let $\phi(t)$ be the unique sharp type traveling wave
corresponding to the speed $c=c^*(m,p,r)$ according to Lemma \ref{le-cstar}.
This kind of special solution is the unique one that
is strictly increasing on $(0,+\infty)$,
$\phi(+\infty)=K$ and $\phi'(+\infty)=0$,
according to the monotone dependence Lemma \ref{le-dependent} and Lemma \ref{le-cstar}.
Detailed discussion can be found in the proof of Lemma 3.10 in \cite{Non20}.

In the proof of Lemma \ref{le-semi-decay},
we formulate the generalized phase plane \eqref{eq-zphase} and \eqref{eq-zphasepsi},
where $\phi_{cr}$ is defined by \eqref{eq-phicr}.
Moreover, $\tilde\psi(0)=0$, $\tilde\psi(K)=0$ since $\phi'(+\infty)=0$,
$\tilde\psi(\phi)>0$ for all $\phi\in(0,K)$.
We rewrite \eqref{eq-zphasepsi} into
\begin{equation} \label{eq-zpsi-vc}
\frac{\mathrm{d}\tilde\psi}{\mathrm{d}\phi}=c
-\frac{m\phi^{m-1}(b(\phi)-d(\phi))}{\tilde\psi^\frac{1}{p-1}}
+\frac{m\phi^{m-1}(b(\phi)-b(\phi_{cr}))}{\tilde\psi^\frac{1}{p-1}},
~\phi\in(0,K).
\end{equation}
For any $g\in\mathscr{D}=\{g\in C^1([0,K]);\int_0^K g(s)\mathrm{d}s=1,g(s)>0,g'(s)<0,\forall s\in(0,K)\}$,
we multiply \eqref{eq-zpsi-vc} by $g(s)$ and integrate it over $(0,K)$ to find
\begin{align} \nonumber
c=&\int_0^K g(\phi)\frac{\mathrm{d}\tilde\psi}{\mathrm{d}\phi}\mathrm{d}\phi
+\int_0^K g(\phi)\frac{m\phi^{m-1}(b(\phi)-d(\phi))}{\tilde\psi^\frac{1}{p-1}}\mathrm{d}\phi
\\ \nonumber
&-\int_0^K g(\phi)\frac{m\phi^{m-1}(b(\phi)-b(\phi_{cr}))}{\tilde\psi^\frac{1}{p-1}}\mathrm{d}\phi
\\ \nonumber
=&\int_0^K -g'(\phi)\tilde\psi(\phi)\mathrm{d}\phi
+\int_0^K g(\phi)\frac{m\phi^{m-1}(b(\phi)-d(\phi))}{\tilde\psi^\frac{1}{p-1}}\mathrm{d}\phi
\\ \nonumber
&+\big[g(\phi)\tilde\psi(\phi)\big]\Big|_{\phi=0}^{\phi=K}
-\int_0^K g(\phi)\frac{m\phi^{m-1}(b(\phi)-b(\phi_{cr}))}{\tilde\psi^\frac{1}{p-1}}\mathrm{d}\phi
\\ \label{eq-zvar}
=&
F(g,\tilde\psi)
-\int_0^K g(\phi)\frac{Dm\phi^{m-1}(b(\phi)-b(\tilde\phi_{cr}(\phi)))}{\tilde\psi}\mathrm{d}\phi,
\end{align}
according to
$\int_0^K g(s)\mathrm{d}s=1$ and $[g(\phi)\tilde\psi(\phi)]\big|_{\phi=0}^{\phi=K}=0$ as
$\tilde\psi(0)=0=\tilde\psi(K)$,
where the functional
$$
F(g,\tilde\psi):=\int_0^K -g'(\phi)\tilde\psi(\phi)\mathrm{d}\phi
+\int_0^K g(\phi)\frac{m\phi^{m-1}(b(\phi)-d(\phi))}{\tilde\psi^\frac{1}{p-1}}\mathrm{d}\phi.
$$

Next, we consider the function $F(g,\tilde\psi)$ over all set $\mathscr{D}$.
Utilizing Young's Inequality, we see that
\begin{equation} \label{eq-functional}
F(g,\tilde\psi)
\ge \int_0^K \frac{p}{(p-1)^{(p-1)/p}}(-g'(\phi))^\frac{1}{p}(g(\phi))^\frac{p-1}{p}
(m\phi^{m-1}(b(\phi)-d(\phi)))^\frac{p-1}{p}\mathrm{d}\phi,
\end{equation}
see also in \cite{Benguria18} for non-delayed case.
The equality in \eqref{eq-functional} is attainable if there exists a function $\hat g\in \mathscr{D}$ such that
\begin{equation} \label{eq-hatg}
m\phi^{m-1}(b(\phi)-d(\phi))\hat g(\phi)=(p-1)(-\hat g'(\phi))\tilde\psi^\frac{p}{p-1}(\phi), \quad \phi\in(0,K).
\end{equation}

The existence of a $\hat g\in\mathscr{D}$ solving \eqref{eq-hatg} relies heavily on the asymptotic behavior
of $\tilde\psi(\phi)$ near $K$ and near $0$
as shown in Lemma \ref{le-asymp} and Lemma \ref{le-semi-1}.
According to Lemma \ref{le-asymp},
$$
\tilde\psi^\frac{p}{p-1}(\phi)
\sim
\begin{cases}
\kappa (K-\phi)^p, \quad & p\ge2,\\
\kappa (K-\phi)^2, \quad & p\in(1,2),
\end{cases}
$$
where $\kappa>0$ is a positive number.
Therefore, $\hat g(K)=0$.
Otherwise,
\begin{equation} \label{eq-zhatg}
\frac{-\hat g'(\phi)}{\hat g(\phi)}=\frac{m\phi^{m-1}(b(\phi)-d(\phi))}{(p-1)\tilde\psi^\frac{p}{p-1}(\phi)}
\sim \frac{mK^{m-1}(d'(K)-b'(K))}{(p-1)\kappa (K-\phi)^{\max\{p-1,1\}}},
\quad \text{as~} \ \phi\to K^-,
\end{equation}
which is not integrable, a contradiction.
Consider the singular ODE \eqref{eq-zhatg} near $K$ with the condition $\hat g(K)=0$, and $\hat g'(\phi)<0$ for $\phi\in(0,K)$,
since $\max\{p-1,1\}\ge1$ for all $p>1$, it has infinitely many solutions such that $\hat g(\phi)>0$ for $\phi\in(0,K)$.
If $\hat g(0)$ is finite, we can normalize $\hat g$ such that $\int_0^K\hat g(\phi)\mathrm{d}\phi=1$ and then $\hat g\in\mathscr{D}$.
According to Lemma \ref{le-semi-1}, $\tilde\psi(\phi)\sim c\phi$ as $\phi\to0^+$.
Then
\begin{equation} \label{eq-zhatg0}
\frac{-\hat g'(\phi)}{\hat g(\phi)}=\frac{m\phi^{m-1}(b(\phi)-d(\phi))}{(p-1)\tilde\psi^\frac{p}{p-1}(\phi)}
\sim \frac{m(b'(0)-d'(0))\phi^{m}}{(p-1)c^\frac{p}{p-1}\phi^\frac{p}{p-1}},
\quad \text{as~} \ \phi\to K^-.
\end{equation}
It follows that $\hat g(0)<+\infty$ since $m-\frac{p}{p-1}>-1$ as
$m(p-1)>1$ such that $\phi^{m-\frac{p}{p-1}}$ is integrable near zero.

Finally, for $\hat g\in\mathscr{D}$, we have
\begin{align*}
c=&F(\hat g,\tilde\psi)
-\int_0^K \hat g(\phi)\frac{Dm\phi^{m-1}(b(\phi)-b(\tilde\phi_{cr}(\phi)))}{\tilde\psi}\mathrm{d}\phi \\
<&F(\hat g,\tilde\psi) \\
\le&\sup_{g\in\mathscr{D}}\int_0^K \frac{p}{(p-1)^{(p-1)/p}}(-g'(\phi))^\frac{1}{p}(g(\phi))^\frac{p-1}{p}
(m\phi^{m-1}(b(\phi)-d(\phi)))^\frac{p-1}{p}\mathrm{d}\phi \\
=&c^*(m,p,0),
\end{align*}
where the last equality is the variational characterization of the speed for non-delayed case
as proved in \cite{Benguria18}.
The proof is completed.
$\hfill\Box$

{\it\bfseries Proof of Theorem \ref{th-existence}.}
The existence and uniqueness of sharp type traveling wave
are proved in Lemma \ref{le-cstar}.
According to Lemma \ref{le-cstarr}, we see that time delay slows down the critical wave speed.
Any other traveling waves must be positive and the regularity is trivial since
\eqref{eq-tw} is non-degenerate at where $\phi(t)>0$.
$\hfill\Box$

{\it\bfseries Proof of Theorem \ref{th-sharp}.}
This is proved according to the asymptotic behavior near $0$ in Lemma \ref{le-semi-1}.
$\hfill\Box$

\section{Model Formulation}

The models with degenerate diffusion but without time-delay were firstly introduced in \cite{Gurney75,Murry},
and the models with time-delay and regular diffusion were widely studied in \cite{Chern,LLLM,Mei_LinJDE09,MLSS,Mei-Ou-Zhao,Schaaf,So-Yang,So-Wu-Zou}. However, the derivation of the models with both effects of degenerate diffusion and time-delay is not officially derived, even if such  models had been proposed and studied in our previous research works \cite{HJMY,JDE18,Non20,JDE20} based on the mathematical concerning.
In this section, we develop a degenerate diffusion model with  time delay that arises in the modelling of age-structured populations. Here, we give a brief derivation of the equations we treat.

The problem is as follows. Let $a$ denote chronological age, $t$ denote time and $x$ denote
spatial position, and let $u(a, t, x)$ denote the population density of age $a$ at time $t$ and at position $x$.
Here,  we concerns species whose life cycles consist of two demographically
distinct phases incorporating immature and reproductive periods.
 By $r\ge 0$ we denote the maturation time that divides the two phases,
so the matured population density at location $x$ and time $t$ is
\begin{equation}\label{eq-P}
u(x,t)=\int_r^\infty v(t,x,a) da.
\end{equation}
Also, since only the mature can reproduce,
the functional dependence of birth rate $\beta$ is assumed to enter only through
dependence on $u$ and so that $\beta=\beta(u)$.

Assuming that the emigration in species due to intraspecific competition in a way that makes
the flux of individuals proportional to the gradient of the mature population, in \cite{Busenberg}, the following
age-structured population model with degenerate diffusion is derived
\begin{align}\label{eq-age-sturcture}
\begin{cases}
\displaystyle
\pd v a+\pd v t=g(t,x,u,\nabla u)\cdot \nabla v+h(t,x,u,\nabla u,D^2 u)v-\mu v,\quad &x\in \mathbb R^n,~t>0,
\\[1mm]
v(0,t,x)=\int_r^\infty \beta(u)v(a,t,x)da,
\\[1mm]
v(a,0,x)=v_0(a,x).
\end{cases}
\end{align}
Using \eqref{eq-P} and integrating the partial differential equation \eqref{eq-age-sturcture}  from $r$ to $+\infty$, we obtain for $t>0$ and $x\in \mathbb R^n$
\begin{equation}\label{eq-delaymodel}
\pd u t=\tilde g(t,x,u,\nabla u)\cdot\nabla u+\tilde h(t,x,u,\nabla u,D^2u)u
+\beta (u(t-r,x))u(t-r,x)S-\mu u.
\end{equation}
In obtaining \eqref{eq-delaymodel} the following biological realistic assumptions are necessary \cite{Sulsky}:
(i) $v(a,t,x)\to 0$ as $a\to +\infty$;
(ii) the birth rate at time $t$ reduces to
$v(0,t,x)=\beta(u) u$;
and (iii) $v(r,t,x)=v(0,t-r,x)S$,
where $S$ is the fraction of individuals that survives through
the first demographic phase.
The last point reflects the fact that the individuals of age $r$ is made up of the survivors that were born at time
$t-r$.

The spatial diffusion term $\tilde g\cdot\nabla u+\tilde hu$ includes operators of the form $\Delta u^m$ and doubly nonlinear operator in \eqref{eq-main}.
When $r=0, p=2, m>1$, \eqref{eq-delaymodel} may reduce to
density dependent models of population dynamics
of the form
\begin{equation} \label{eq-densitynondelay}
\pd u t=\Delta u^m+\beta(u)u-\mu u,\quad x \in \mathbb R^n,~t>0,
\end{equation}
which was considered by means of existing theory (Aronson \cite{Aronson80Density},
Gurtin \& MacCamy \cite{Gurtin77}).
Our model includes a large number of evolution equations in biology, for example,
the degenerate delayed Fisher-KPP equations,
the degenerate Nicholson's blowflies equation,
and the Mackey-Glass equation.

\section*{Acknowledgement}

The research of S. Ji was supported
by Guangdong Basic and Applied Basic Research Foundation
and the Fundamental Research
Funds for the Central Universities of SCUT.
The research of M. Mei was supported in part by
NSERC Grant RGPIN 354724-16, and FRQNT Grant No. 2019-CO-256440.
The research of T. Xu was supported in part by NSFC Grant No. 11871230.
The research of J. Yin was supported in part by NSFC Grant No. 11771156,
NSF of Guangzhou Grant No. 201804010391,
and Guangdong Basic and Applied Basic Research Foundation Grant No.2020B1515310013.


\begin{thebibliography}{99}






\bibitem{Aronson80Density}
D.G. Aronson,
Density-dependent interaction-diffusion systems, in: Proc. Adv. Seminar on Dynamics and Modeling
of Reactive System, Academic Press, New York, 1980.


\bibitem{VazquezJDE17}
A. Audrito and J.L. V\'azquez,
The Fisher-KPP problem with doubly nonlinear diffusion,
\newblock {\em J. Differential Equations},
{\bf 263} (2017), 7647--7708.

\bibitem{Benguria-Variational}
R.D. Benguria and M.C. Depassier,
Variational principle for the asymptotic speed of fronts of the density-dependent diffusion-reaction equation,
{\em Physical Review E},
{\bf 52} (1995), 3285--3287.


\bibitem{Benguria}
R.D. Benguria and M.C. Depassier,
Variational characterization of the speed of propagation of fronts for the nonlinear diffusion equation,
\newblock {\em Commun. Math. Phys.},
{\bf 175} (1996), 221--227.

\bibitem{Benguria18}
R.D. Benguria, M.C. Depassier,
Variational characterization of the speed of reaction diffusion fronts for gradient dependent diffusion,
{\em Ann. Henri Poincar\'e}, {\bf 19} (2018), 2717--2726.


\bibitem{Busenberg}
S. Busenberg, M. Iannelli,
A class of nonlinear diffusion problems in age-dependent population dynamics,
{\it Nonlinear Analysis Theory Methods \& Applications}, {\bf 7} (1983), 501--529.

\bibitem{Chern}
I.-L. Chern, M. Mei, X. Yang and Q. Zhang,
Stability of non-monotone critical traveling waves for reaction-diffusion equations with time-delay,
{\em J. Differential Equations},
{\bf 259} (2015), 1503--1541.

\bibitem{Faria}
T. Faria and S. Trofimchuk,
Non-monotone travelling waves in a single species reaction-diffusion equation with delay,
{\em J. Differential Equations},
{\bf 228} (2006), 357--376.

\bibitem{Fisher}
R.A. Fisher,
The wave of advance of advantageous genes,
{\em Ann. Eugen.},
{\bf 7} (1937), 335--369.



\bibitem{Gomez}
A. Gomez and S. Trofimchuk,
Global continuation of monotone wavefronts,
{\em J. Lond. Math. Soc.},
{\bf 89} (2014), 47--68.

\bibitem{Gurney80}
W.S.C. Gurney, S.P. Blythe and R.M. Nisbet,
Nicholson's blowflies revisited,
{\em Nature},
{\bf 287} (1980), 17--21.


\bibitem{Gurney75}
W.S.C. Gurney and R.M. Nisbet,
The regulation of inhomogeneous population,
{\em J. Theors. Biol.},
{\bf 52} (1975), 441--457.


\bibitem{Gurtin77}
M.E. Gurtin, R.C. MacCamy, On the diffusion of biological populations,
{\it Mathematical Biosciences}, {\bf 33} (1977), 35--49.


\bibitem{HJMY}
R. Huang, C.H. Jin, M. Mei and J.X. Yin,
Existence and stability of traveling waves for degenerate reaction-diffusion equation with time delay,
{\em J. Nonlinear Sci.},
{\bf 28} (2018), 1011--1042.


\bibitem{Jin-JDE}
C.H. Jin, T. Yang, J.X. Yin,
Waiting time for a non-Newtonian polytropic filtration equation with convection,
{\it J. Differential Equations}, {\bf 252} (2012), 4862--4885



\bibitem{LLLM}
C.-K. Lin, C.-T. Lin, Y. Lin, M. Mei,
Exponential stability of nonmonotone traveling waves for Nicholson's
blowflies equation, {\it  SIAM J. Math. Anal.}, {\bf 46} (2014), 1053--1084.



\bibitem{Mackey}
M.C. Mackey and L. Glass,
Oscillation and chaos in physiological control systems,
{\em Science},
{\bf 197} (1977), 287--289.

\bibitem{Mei_LinJDE09}
M. Mei, C.K. Lin, C.T. Lin and J.W.-H. So,
Traveling wavefronts for time-delayed reaction-diffusion equation: (I) local nonlinearity,
{\em J. Differential Equations},
{\bf 247} (2009), 495--510.


\bibitem{MLSS}
M. Mei, J. W.-H. So, M.Y. Li, S. S.P. Shen,
Asymptotic stability of travelling waves for Nicholson's blowflies
equation with diffusion, {\it Proc. R. Soc. Edinb.}, {\bf 134A} (2004), 579--594.

\bibitem{Mei-Ou-Zhao}
M. Mei, C. Ou, X.-Q. Zhao,
Global stability of monostable traveling waves for nonlocal time-delayed reaction-diffusion equations,
{\it SIAM J. Math. Anal.},  {\bf 42} (2010), 2762--2790.



\bibitem{Murry}
J.D. Murry,
Mathematical Biology I: An Introduction,
Springer, New York, USA,
2002.




\bibitem{Schaaf}
K.W. Schaaf,
Asymptotic behavior and traveling wave solutions for parabolic functional differential equations,
{\em Trans. Amer. Math. Soc.},
{\bf 302} (1987), 587--615.

\bibitem{So-Wu-Zou}
J.W.-H. So, J. Wu, X. Zou,
A reaction-diffusion model for a single species with age
structure: (I) Traveling wavefronts on unbounded domains,
{\it Proc. R. Soc. Lond. Ser. A Math. Phys. Eng. Sci.}, {\bf 457} (2001),  1841--1853.


\bibitem{So-Yang}
J.W.-H. So, Y. Yang,
Dirichlet Problem for the diffusion Nicholson's blowflies equation,
{\it J. Differential Equations}, {\bf 150} (1998), 317--348.

\bibitem{Sulsky}
D. Sulsky, R.R. Vance, W.I. Newman,
Time delays in age-structured populations,
{\it J. Theor. Biol.}, {\bf 141} (1990), 403--422.



\bibitem{JDE18}
T.Y. Xu, S.M. Ji, M. Mei and J.X. Yin,
Traveling waves for time-delayed reaction diffusion equations with degenerate diffusion,
{\em J. Differential Equations},
{\bf 265} (2018), 4442--4485.

\bibitem{Non20}
T.Y. Xu, S.M. Ji, M. Mei, and J.X. Yin,
Variational approach of critical sharp front speeds in degenerate diffusion model with time delay,
{\em Nonlinearity}, {\bf 33} (2020), 4013--4029.

\bibitem{JDE20}
T.Y. Xu, S.M. Ji, M. Mei, and J.X. Yin,
Sharp oscillatory traveling waves of structured population dynamics model with degenerate diffusion,
{\em J. Differential Equations}, {\bf 269} (2020), 8882--8917.

\bibitem{ArXiv21}
T.Y. Xu, S.M. Ji, M. Mei, and J.X. Yin,
Propagation speed of degenerate diffusion equations
with time delay, preprint 2021. ArXiv:2011.14813

\end{thebibliography}
\end{document}